\documentclass[11pt,reqno]{amsart}
 \usepackage[cp1251]{inputenc}
\usepackage[english]{babel}
\usepackage{amsmath}
\usepackage{amssymb}
\usepackage{amsfonts}
\usepackage{xcolor}
\usepackage[colorlinks]{hyperref}
 \usepackage{amsmath,systeme}
\usepackage{nicefrac}
\usepackage{latexsym,fancyhdr,amssymb,color,amsmath,amsthm,graphicx,listings,comment}

\usepackage{titlesec}
\titlespacing{\section}{0pt}{2.5ex}{1.5ex}
\titlespacing{\subsection}{0pt}{1.5ex}{1ex}
\titlespacing{\subsubsection}{0pt}{1.5ex}{1ex}
\titleformat{\section}{\large\bfseries\centering}{\thesection}{1em}{}
\titleformat{\subsection}[runin]{\bfseries}{\thesubsection.}{0.5em}{}[.\mbox{\ }]
\titleformat{\subsubsection}[runin]{\bfseries}{\thesubsubsection.}{0.4em}{}[.\mbox{\ }]

\newtheorem{conjecture}{Conjecture}

\newtheorem{theorem}{Theorem}

\newtheorem{corollary}{Corollary}

\newtheorem*{example}{Examples}
\newtheorem{problem}{Problem}

\definecolor{sangria}{rgb}{0.57, 0.0, 0.04}

\hypersetup{linkcolor=red,filecolor=black,urlcolor=blue, citecolor=blue}

\begin{document}
\renewcommand{\refname}{References}
\renewcommand{\proofname}{Proof.}
\renewcommand{\figurename}{Fig.}

\title[ ]{ \large Primitive Euler brick generator}

  \author[Djamel HIMANE]{{Djamel Himane}   }

\address{\textbf{Djamel  Himane}   
\newline\hphantom{iii} LA3C Laboratory, Faculty of Mathematics,
\newline\hphantom{iii}  USTHB,  Algiers, Algeria.}%
\email{\textcolor{blue}{dhimane@usthb.dz}}%

 \maketitle { \small 
 

\begin{quote}
\noindent{\bf Abstract:}  The smallest Euler brick, discovered by Paul Halcke, has edges
 $(177, 44, 240) $ and face diagonals $(125, 267, 244 ) $, generated by the primitive Pythagorean triple $ (3, 4, 5). $  
 Let $ (u,v,w) $ primitive Pythagorean triple, Sounderson made a generalization parameterization of the edges
\begin{equation*}  
 a  =   \vert u(4v^2 - w^2) \vert, \quad b = \vert v(4u^2 - w^2)\vert, \quad c = \vert 4uvw \vert  
\end{equation*}
give face diagonals 
\begin{equation*}
 {\displaystyle d=w^{3},\quad e=u(4v^{2}+w^{2}),\quad f=v(4u^{2}+w^{2})}
\end{equation*}
 leads to an Euler brick. Finding other formulas that generate these primitive bricks, other than formula above, or making initial guesses that can be improved later, is the key to understanding how they are generated.

 \end{quote} }
 
   \vspace{.5cm}

\section*{Introduction} 

An Euler brick is a cuboid of integer side dimensions $ a, b, c $ such that the face diagonals are
integers.  In geometric terms is equivalent to a solution to the following system of Diophantine equations:
\begin{equation} \label{eq01} 
{\displaystyle {\begin{cases}
\hspace{1cm} a^{2}+b^{2}=d^{2} \\
\hspace{1cm} a^{2}+c^{2}=e^{2}\\
\hspace{1cm} b^{2}+c^{2}=f^{2}
\end{cases}}}
\end{equation}
where $ a, b, c $ are the edges and $ d, e, f $ are the diagonals.

 Given an Euler brick with edge-lengths $ (a, b, c) $, the triple $ (bc, ac, ab) $ constitutes an Euler brick as well \cite{4}.
 If $ \gcd(a, b, c)=1 $, we call it a primitive Euler  brick.

A perfect Euler brick  is an Euler brick whose space diagonal also has integer length. In other words, the following equation is added to the system of Diophantine equations defining an Euler brick:
\begin{equation}
{\displaystyle a^{2}+b^{2}+c^{2}=g^{2},}
\end{equation}
where $ g $ is the space diagonal. As of May 2024 \cite{5}, no example of a perfect cuboid had been found and no one has proven that none exist.
\vspace{0.5cm}

\parbox{14cm}{
\scalebox{0.09}{ \includegraphics{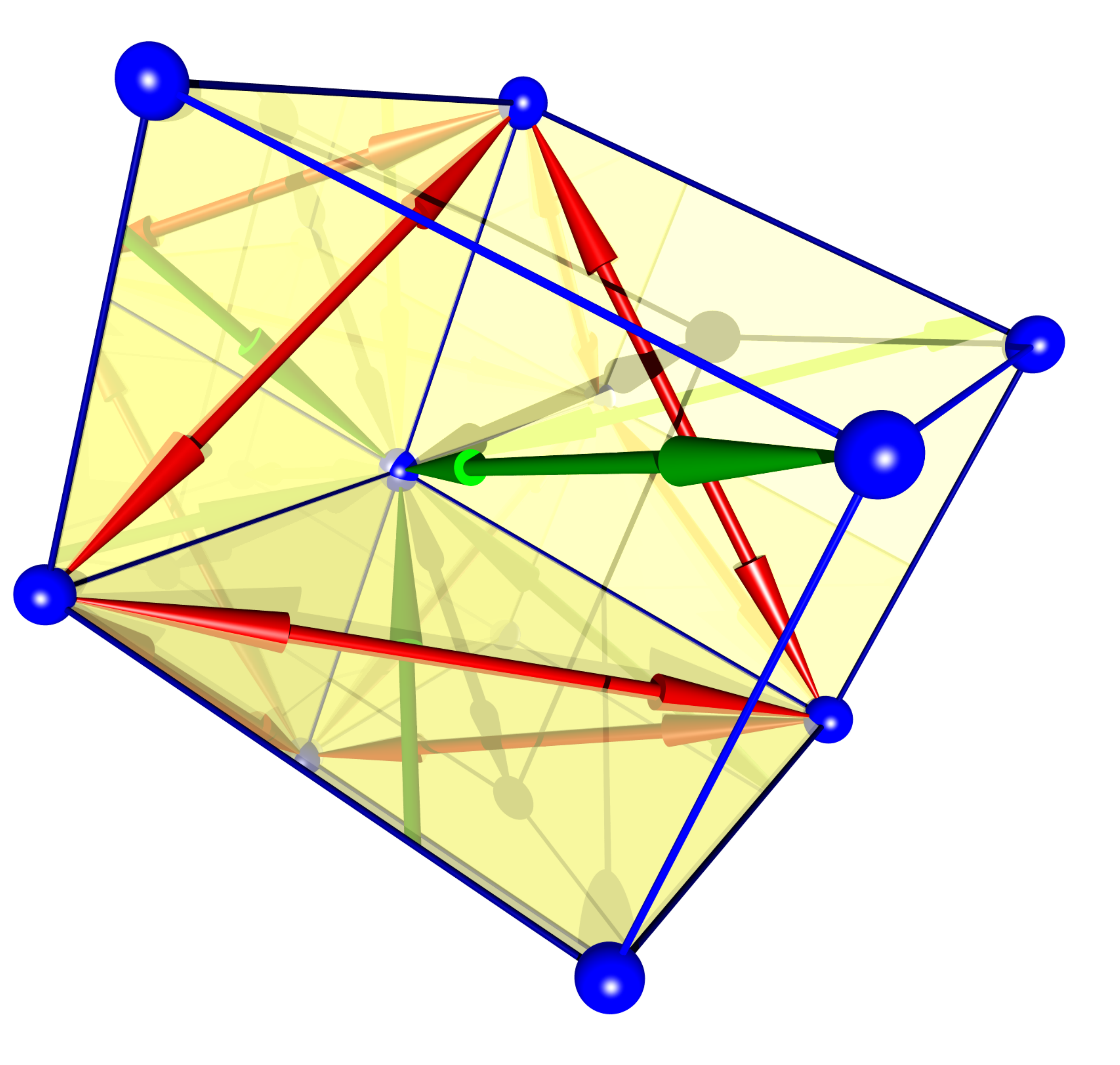}}
\hspace{1cm}
\scalebox{0.55}{ \includegraphics{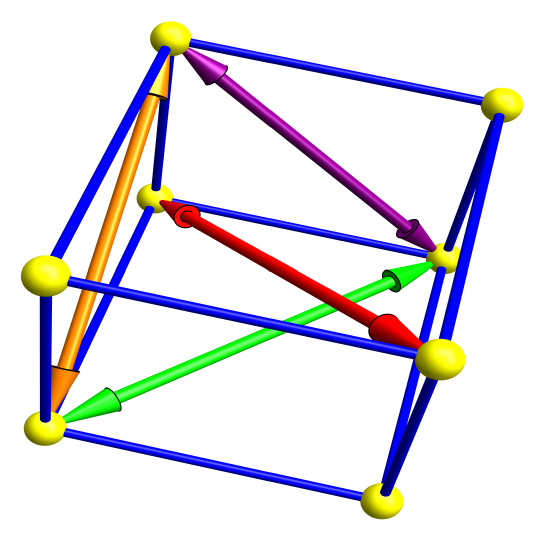}}
}

\begin{small}
\begin{center}
 The figures included are by  Oliver Knill, showing 3D drawings of an Euler brick \cite{2}  
\end{center}
 \end{small}
  \vspace{0.5cm}

 We emphasize that in everything that follows, sides $ a, b, c $ are arranged according to the powers of $  2$ that they contain, that is to say 
\begin{equation}
0= v_{2}(a)\quad < \quad v_{2}(b) \quad  < \quad  v_{2}(c)
\end{equation}
If $ u, v, w $ are integers satisfying $ u^2 + v^2 = w^2 $, then the Sounderson \cite{1} parametrization of the edges
\begin{equation} \label{eq02} 
 a  =   \vert u(4v^2 - w^2) \vert, \quad b = \vert v(4u^2 - w^2)\vert, \quad c = \vert 4uvw \vert  
\end{equation}
give face diagonals
\begin{equation} \label{eq03} 
{\displaystyle d=w^{3},\quad e=u(4v^{2}+w^{2}),\quad f=v(4u^{2}+w^{2}).}
\end{equation}
leads to an Euler brick.
The smallest Euler brick, discovered by Paul Halcke in 1719 \cite{3,4}, has edges $ (a, b, c) = (177, 44, 240) $ and face diagonals $ (d, e, f ) = (125, 267, 244 ) $, generated by the primitive Pythagorean triple $ (3, 4, 5). $

There are many Euler bricks which are not parametrized as above, for instance the Euler brick with edges $ (a, b, c) = (275, 252, 240) $ and face diagonals $ (d, e, f ) = (373, 365, 348). $

Finding other formulas that generate these primitive bricks, other than formula \eqref{eq02}, \eqref{eq03}, or making initial guesses that can be improved later, is the key to understanding how they are generated.
 
Euclid's formula says that, $(u,v,w)$ are a primitive Pythagorean triple,i.e.  
\begin{equation*}
 u^2 + v^2=w^2
 \end{equation*} for $ u,v,w $ are integers, if and only if \cite{4}
 \begin{equation*}
 u=m^2-n^2, \hspace{2cm} v=2mn,  \hspace{2cm} w=m^2 + n^2
\end{equation*}   for some integers $ m,n $ satisfait 
  \begin{equation*}
  m > n > 0, \hspace{1cm} \gcd(m, n) = 1, \hspace{1cm}  m + n \equiv 1 \pmod 2 .
 \end{equation*}
 Consider the following the primitive Pythagorean triples $ (u_{i}, v_{i}, w_{i}) $,  
 \begin{equation} \label{eq04}
  \text{for} \;   i=0,1,2,  \quad   \text{with} \;  u_{i} \; \text{odd and} \;  v_{i} \; \text{ even},
   \quad \text{where}  \quad  u_{0}^{2}+v_{0}^{2}= t_{0} ,
 \end{equation}
    with $ u_{0}$ odd, $ v_{0} $ it is not necessarily even and $ t_{0} $  it is not necessarily a square.

 \section*{Main results}

\begin{theorem} \label{th1} 
 Let $ (u_{1}, v_{1}, w_{1}) $ and $ (u_{2}, v_{2}, w_{2}) $ two primitive Pythagorean triples and $ (u_{0}, v_{0}) $  two positive integers satisfy the condition \eqref{eq04}, where $ (u_{0}u_{2})^2+(v_{0}v_{1})^2=e^2=\square $ and $ v_{0}u_{1} = u_{0} v_{2} $, then the parametrization of the edges
\begin{equation} \label{eq05} 
 a  =    u_{0}u_{2} , \hspace{2cm} b = v_{0}u_{1} = u_{0} v_{2}, \hspace{2cm} c = v_{0} v_{1}
\end{equation}
give face diagonals
\begin{equation} \label{eq06} 
{\displaystyle d=u_{0}w_{2},\hspace{1cm} e= \sqrt{(u_{0}u_{2})^2+(v_{0}v_{1})^2 } ,\hspace{1cm} f = v_{0}w_{1}  }
\end{equation}
leads to an Euler brick. 
 \end{theorem}
 
 \begin{proof}
 It is clear that if $ (u_{1}, v_{1}, w_{1}) $ , $ (u_{2}, v_{2}, w_{2}) $ and  $ (u_{0}, v_{0}) $  satisfy the conditions mentioned above, then the sides $ a, b, c $ and the diagonal faces $ d, e, f $ are positive integers.
 All that remains is to prove that, it is a solution to the system of equations \eqref{eq01}, we have:
 \begin{equation*}  
{\displaystyle {\begin{cases}
\hspace{.5cm} a^{2}+b^{2} = (u_{0}u_{2})^2 + (u_{0}v_{2})^2 = u_{0}^{2}(u_{2}^2 + v_{2}^2)=(u_{0}w_{2})^2=d^{2} \\
\hspace{.5cm} a^{2}+c^{2}=(u_{0}u_{2})^2+(v_{0}v_{1})^2=e^2 \hspace{2cm} \text{(hypothesis)} \\
\hspace{.5cm} b^{2}+c^{2}=(v_{0}u_{1})^2 + (v_{0}v_{1})^2 = v_{0}^{2}(u_{1}^2 + v_{1}^2)=(v_{0}w_{1})^2=f^{2}
\end{cases}}}
\end{equation*}
 \end{proof}
 
 \begin{table}[h!]
  \begin{tabular}{|c||c|c||c|c|c||c|c|c||c|c|c||c|c|c| } 
\hline
$ \# $ & $ u_{0} $ & $ v_{0} $ & $ u_{1} $ & $ v_{1} $ & $ w_{1} $ & $ u_{2} $ & $ v_{2} $ & $ w_{2} $ & $ a $ & $ b $ & $ c$  & $ d $ & $e $ & $ f$  \\ \hline
1  & 7 & 20 & 7  &  24  &  25 &  99 & 20 & 101 &  693 &  140 &  480 &  707 &  843 &  500     \\ \hline
2  & 33 & 8 & 99 &  20 & 101 & 7 & 24 & 25 & 231 & 792&160& 825 & 281 &   808  \\ \hline
3 & 23 & 12 & 483 & 44 & 485 & 275  & 252  &  373 & 6325 & 5796 & 528 & 8579 &6347&  5820    \\ \hline
4  & 17 & 12 & 85 & 132 & 157 & 11 & 60 & 61 & 187 & 1020 &1584 & 1037  & 1595 &    1884     \\ \hline
5  & 1 & 4 & 11 & 60 & 61 & 117 & 44 & 125 & 117 & 44 & 240 & 125 & 267  & 244 \\ \hline
 
\end{tabular}
\bigskip
  \caption{The first five primitive Euler brick, given by the formula \eqref{eq05}, \eqref{eq06} } 
 \end{table}

From the formula \eqref{eq05}, \eqref{eq06} concludes $ u_{0} \mid u_{1} $ , $ v_{0}\mid v_{2} $ and deduce that

 \begin{corollary}
 Let $ (u_{1}, v_{1}, w_{1}) $ and $ (u_{2}, v_{2}, w_{2}) $ two primitive Pythagorean triples satisfy the condition \eqref{eq04}, where $ (u_{1}u_{2})^2+(v_{1}v_{2})^2=e^2=\square $, then the parametrization of the edges
\begin{equation} \label{eq07} 
 a  =    u_{1}u_{2} , \hspace{2cm} b =  u_{1} v_{2}, \hspace{2cm} c = v_{1} v_{2}
\end{equation}
give face diagonals
\begin{equation} \label{eq08} 
{\displaystyle d=u_{1}w_{2},\hspace{1cm} e= \sqrt{(u_{1}u_{2})^2+(v_{1}v_{2})^2 } ,\hspace{1cm} f = v_{2}w_{1}  }
\end{equation}
leads to an Euler brick. 
 \end{corollary}
 
 \begin{proof}
 It is clear that if $ (u_{1}, v_{1}, w_{1}) $ and $ (u_{2}, v_{2}, w_{2}) $ satisfy the conditions mentioned above, then the sides $ a, b, c $ and the diagonal faces $ d, e, f $ are positive integers.
 All that remains is to prove that, it is a solution to the system of equations \eqref{eq01}, we have:
 \begin{equation*}  
{\displaystyle {\begin{cases}
\hspace{.5cm} a^{2}+b^{2} = (u_{1}u_{2})^2 + (u_{1}v_{2})^2 = u_{1}^{2}(u_{2}^2 + v_{2}^2)=(u_{1}w_{2})^2=d^{2} \\
\hspace{.5cm} a^{2}+c^{2}=(u_{1}u_{2})^2+(v_{1}v_{2})^2=e^2 \hspace{2cm} \text{(hypothesis)} \\
\hspace{.5cm} b^{2}+c^{2}=(u_{1}v_{2})^2 + (v_{1}v_{2})^2 = v_{2}^{2}(u_{1}^2 + v_{1}^2)=(v_{2}w_{1})^2=f^{2}
\end{cases}}}
\end{equation*}
 \end{proof}
 
 \begin{table}[h!]
  \begin{tabular}{|c||c|c|c||c|c|c||c|c|c||c|c|c| } 
\hline
$ \# $ & $ u_{1} $ & $ v_{1} $ & $ w_{1} $ & $ u_{2} $ & $ v_{2} $ & $ w_{2} $ & $ a $ & $ b $ & $ c$  & $ d $ & $e $ & $ f$  \\ \hline
1  &  7  &  24  &  25 &  99 & 20 & 101 &  693 &  140 &  480 &  707 &  843 &  500            \\ \hline
2          &             &             &             &             &  & &&&&   &&           \\ \hline
3          &           &            &           &             &   &&&&&     &&      \\ \hline
4          &            &             &       &      &    &&&&&   &&         \\ \hline
5          &           &          &         &       &  &&&&&&  &    \\ \hline
 
\end{tabular}
\bigskip
  \caption{The first primitive Euler brick, given by the formula \eqref{eq07}, \eqref{eq08} } 
 \end{table}
 
\begin{problem}
Let $ (u_{1}, v_{1}, w_{1}) $ and $ (u_{2}, v_{2}, w_{2}) $ two primitive Pythagorean triples and $ (u_{0}, v_{0}) $  two positive integers, related to Theorem \ref{th1}, If the condition $ (u_{0}u_{2})^2+(v_{0}v_{1})^2=\square $ is met, then $ v_{0}u_{1} = u_{0} v_{2} $ ?    
\end{problem}

\begin{theorem} \label{th2} 
 Let $ (u_{1}, v_{1}, w_{1}) $ and $ (u_{2}, v_{2}, w_{2}) $ two primitive Pythagorean triples and $ (u_{0}, v_{0}) $  two odd positive integers satisfy the condition \eqref{eq04}, where $ (v_{0}v_{1})^2+(u_{0}v_{2})^2=f^2=\square $ and $ u_{0}u_{2} = v_{0}u_{1} $, then the parametrization of the edges
\begin{equation} \label{eq09} 
 a  =  u_{0}u_{2} = v_{0}u_{1} , \hspace{2cm} b = v_{0} v_{1}, \hspace{2cm} c = u_{0} v_{2} \hspace{1cm}
\end{equation}
give face diagonals
\begin{equation} \label{eq10} 
{\displaystyle d = v_{0}w_{1},\hspace{1cm} e = u_{0}w_{2} ,\hspace{1cm} f = \sqrt{(v_{0}v_{1})^2+(u_{0}v_{2})^2 }  }
\end{equation}
leads to an Euler brick. 
 \end{theorem}
 
 \begin{proof}
 It is clear that if $ (u_{1}, v_{1}, w_{1}) $ , $ (u_{2}, v_{2}, w_{2}) $ and  $ (u_{0}, v_{0}) $  satisfy the conditions mentioned above, then the sides $ a, b, c $ and the diagonal faces $ d, e, f $ are positive integers.
 All that remains is to prove that, it is a solution to the system of equations \eqref{eq01}, we have:
 \begin{equation*}  
{\displaystyle {\begin{cases}
\hspace{.5cm} a^{2}+b^{2} = (v_{0}u_{1})^2 + (v_{0}v_{1})^2 = v_{0}^{2}(u_{1}^2 + v_{1}^2)=(v_{0}w_{1})^2=d^{2} \\
\hspace{.5cm} a^{2}+c^{2}= (u_{0}u_{2})^2 + (u_{0}v_{2})^2 = u_{0}^{2}(u_{2}^2 + v_{2}^2)=(u_{0}w_{2})^2=e^{2} \\
\hspace{.5cm} b^{2}+c^{2}=(v_{0}v_{1})^2+(u_{0}v_{2})^2=f^2 \hspace{2cm} \text{(hypothesis)}
\end{cases}}}
\end{equation*}
 \end{proof}

 \begin{table}[h!]
  \begin{tabular}{|c||c|c||c|c|c||c|c|c||c|c|c||c|c|c| } 
\hline
$ \# $ & $ u_{0} $ & $ v_{0} $ & $ u_{1} $ & $ v_{1} $ & $ w_{1} $ & $ u_{2} $ & $ v_{2} $ & $ w_{2} $ & $ a $ & $ b $ & $ c$  & $ d $ & $e $ & $ f$  \\ \hline
1  & 11 & 17 & 11  &  60  &  61 &  17 & 144 & 145 &  187 &  1020 & 1584 & 1037 & 1595 &  1884   \\ \hline
2  & 21  &  55 &  21  &  20  & 29 & 55 & 48 & 73 & 1155 & 1100 & 1008 & 1595 & 1533 &  1492   \\ \hline
3 & 11 & 39 & 11 & 60 & 61 & 39 &80& 89 & 429 & 2340 & 880 &  2379   & 979 & 2500 \\ \hline
4  &  11 & 23 & 275 & 252& 373 & 575 & 48&577& 6347 & 5796 & 528  & 8579 & 6347   & 5820     \\ \hline
5 & 21 & 25 & 483 & 44 & 485 & 575 & 48 & 577 & 12075 & 1100 & 1008 & 12125 & 12117 &  1492  \\ \hline
 
\end{tabular}
\bigskip
  \caption{The first five primitive Euler brick, given by the formula \eqref{eq09}, \eqref{eq10} } 
 \end{table}

From the formula \eqref{eq09}, \eqref{eq10} concludes $ u_{0} \mid u_{1} $ , $ v_{0}\mid u_{2} $ and deduce that

 \begin{corollary}
 Let $ (u_{1}, v_{1}, w_{1}) $ and $ (u_{2}, v_{2}, w_{2}) $ two primitive Pythagorean triples satisfy the condition \eqref{eq04}, where $ (u_{2}v_{1})^2+(u_{1}v_{2})^2=f^2=\square $, then the parametrization of the edges
\begin{equation} \label{eq11} 
 a  =    u_{1}u_{2} , \hspace{2cm} b =  v_{1} u_{2} , \hspace{2cm} c = u_{1} v_{2}
\end{equation}
give face diagonals
\begin{equation} \label{eq12}  
{\displaystyle d=u_{2}w_{1},\hspace{1cm} e=u_{1}w_{2},\hspace{1cm} f =  \sqrt{(v_{1}u_{2})^2+(u_{1}v_{2})^2 }  }
\end{equation}
leads to an Euler brick. 
 \end{corollary}
 
 \begin{proof}
 It is clear that if $ (u_{1}, v_{1}, w_{1}) $ and $ (u_{2}, v_{2}, w_{2}) $ satisfy the conditions mentioned above, then the sides $ a, b, c $ and the diagonal faces $ d, e, f $ are positive integers.
 All that remains is to prove that, it is a solution to the system of equations \eqref{eq01}, we have:
 \begin{equation*}  
{\displaystyle {\begin{cases}
\hspace{.5cm} a^{2}+b^{2} = (u_{1}u_{2})^2 + (v_{1}u_{2})^2 = u_{2}^{2}(u_{1}^2 + v_{1}^2)=(u_{2}w_{1})^2=d^{2} \\
\hspace{.5cm} a^{2}+c^{2}=(u_{1}u_{2})^2+(u_{1}v_{2})^2= v_{1}^{2}(u_{1}^2 + v_{1}^2)=(v_{1}w_{2})^2=e^2  \\
\hspace{.5cm} b^{2}+c^{2}=(v_{1}u_{2})^2 + (u_{1}v_{2})^2 =f^{2} \hspace{2cm} \text{(hypothesis)}
\end{cases}}}
\end{equation*}
 \end{proof}
 
 \begin{table}[h!]
  \begin{tabular}{|c||c|c|c||c|c|c||c|c|c||c|c|c| } 
\hline
$ \# $ & $ u_{1} $ & $ v_{1} $ & $ w_{1} $ & $ u_{2} $ & $ v_{2} $ & $ w_{2} $ & $ a $ & $ b $ & $ c$  & $ d $ & $e $ & $ f$  \\ \hline
1  &  11  &  60  &  61 &  17 & 144 & 145 &  187 &  1020 &  1584 &  1037 &  1595 &  1884            \\ \hline
2  & 21 &  20  & 29 & 55 & 48 & 73 &1155 & 1100 & 1008 & 1595 & 1533 &  1492  \\ \hline
3 & 11 & 60 &  61 & 39 & 80 & 89 & 429 & 2340 & 880 & 2379    & 979 &  2500    \\ \hline
4 &  275  &  252 & 373 &  575    &  48  & 577 & 158125 & 144900 & 13200 & 214475  & 158675 &  145500       \\ \hline
5  & 483  & 44 & 485 & 575 & 48 & 577 & 277725  & 25300 & 23184 &278875& 278691 &  34316  \\ \hline
 
\end{tabular}
\bigskip
  \caption{The first five Euler brick, given by the formula \eqref{eq11}, \eqref{eq12} } 
 \end{table}
  
 \begin{problem}
Let $ (u_{1}, v_{1}, w_{1}) $ and $ (u_{2}, v_{2}, w_{2}) $ two primitive Pythagorean triples and $ (u_{0}, v_{0}) $  two positive integers, related to Theorem \ref{th2}, If the condition $ (v_{0}v_{1})^2+(u_{0}v_{2})^2=\square $ is met, then $ u_{0}u_{2} = v_{0} u_{1} $ ?    
\end{problem}

 \begin{theorem} \label{th3} 
 Let $ (u_{1}, v_{1}, w_{1}) $ and $ (u_{2}, v_{2}, w_{2}) $ two primitive Pythagorean triples and $ (u_{0}, v_{0}) $  two positive integers satisfy the condition \eqref{eq04}, where $ (u_{0}u_{2})^2 + (v_{0}u_{1})^2 = d^2 = \square $  and $ v_{0} v_{1} = u_{0}v_{2} $, then the parametrization of the edges
\begin{equation} \label{eq13} 
 a  =    u_{0}u_{2} , \hspace{2cm} b =  v_{0} u_{1}, \hspace{2cm} c = v_{0} v_{1} = u_{0}v_{2}
\end{equation}
give face diagonals
\begin{equation} \label{eq14} 
{\displaystyle d= \sqrt{(u_{0}u_{2})^2+(v_{0}u_{1})^2 } ,\hspace{1cm} e = u_{0}w_{2} ,\hspace{1cm} f = v_{0}w_{1}  }
\end{equation}
leads to an Euler brick. 
 \end{theorem}
 
 \begin{proof}
 It is clear that if $ (u_{1}, v_{1}, w_{1}) $, $ (u_{2}, v_{2}, w_{2}) $ and $ (u_{0}, v_{0}) $ satisfy the conditions mentioned above, then the sides $ a, b, c $ and the diagonal faces $ d, e, f $ are positive integers.
 All that remains is to prove that, it is a solution to the system of equations \eqref{eq01}, we have:
 \begin{equation*}  
{\displaystyle {\begin{cases}
\hspace{.5cm} a^{2}+b^{2} = (u_{0}u_{2})^2+(v_{0}u_{1})^2=d^2 \hspace{2cm} \text{(hypothesis)}\\
\hspace{.5cm} a^{2}+c^{2}= (u_{0}u_{2})^2 + (u_{0}v_{2})^2 = u_{0}^{2}(u_{2}^2 + v_{2}^2)=(u_{0}w_{2})^2=e^{2} \\
 \hspace{.5cm} b^{2}+c^{2}=(v_{0}u_{1})^2 + (v_{0}v_{1})^2 = v_{0}^{2}(u_{1}^2 + v_{2}^2)=(v_{0}w_{1})^2=f^{2}
\end{cases}}}
\end{equation*}
 \end{proof}
 
 \begin{table}[h!]
  \begin{tabular}{|c||c|c||c|c|c||c|c|c||c|c|c||c|c|c| } 
\hline
$ \# $ & $ u_{0} $ & $ v_{0} $   & $ u_{1} $ & $ v_{1} $ & $ w_{1} $  & $ u_{2} $ & $ v_{2} $ & $ w_{2} $ & $ a $ & $ b $ & $ c$  & $ d $ & $e $ & $ f$  \\ \hline
1  &  3  &  4 &  11 & 60 & 61 &  39 &  80 &  89 &  117 &  44 &  240 & 125 & 267 &  244     \\ \hline
2  & 5  & 12  &   21 & 20 & 29 & 55 & 48 & 73 & 275  & 252 & 240 & 373 & 365 &  348  \\ \hline
3 & 5  & 12 & 11  & 60 & 61 & 17 & 144 & 145 & 85  & 132  & 720 & 157 & 725 &  732  \\ \hline
4 & 3 & 44 &         17 & 144 & 145 & 65 & 2112 & 2113 &  195 & 748 & 6336  & 773 & 6339 &  6380   \\ \hline
5 & 11 & 20 & 117 & 44 & 125 & 39 & 80 & 89 & 429 & 2340 & 880 & 2379  & 979 & 2500 \\ \hline
 
\end{tabular}
\bigskip
  \caption{The first five primitive Euler brick, given by the formula \eqref{eq13}, \eqref{eq14} } 
 \end{table}

 From the formula \eqref{eq13}, \eqref{eq14} concludes $ u_{0} \mid v_{1} $ , $ v_{0}\mid v_{2} $ and deduce that

 \begin{corollary} [\textbf{Sounderson parametrization}] 
 Let $ (u_{0}, v_{0}, w_{0}) $  primitive Pythagorean triple, then the edges
\begin{equation}  \label{eq17}  
 a  =   \vert u_{0}(4v_{0}^2 - w_{0}^2) \vert, \quad\quad b = \vert v_{0}(4u_{0}^2 - w_{0}^2)\vert, \quad\quad c = \vert 4u_{0}v_{0}w_{0} \vert  
\end{equation}
give face diagonals
\begin{equation}  \label{eq18}
{\displaystyle d=w_{0}^{3},\quad\quad e=u_{0}(4v_{0}^{2}+w_{0}^{2}),\quad\quad f=v_{0}(4u_{0}^{2}+w_{0}^{2})}
\end{equation}
leads to an Euler brick. 
 \end{corollary}
 
 \begin{proof}
 Let $ (u_{0}, v_{0}, w_{0}) $ arbitrary primitive Pythagorean triple, we take
  \begin{equation*} 
  (u_{1}, v_{1}, w_{1}) =( |(2u_{0})^2 - w_{0}^2|,\quad 4 u_{0}w_{0},\quad (2u_{0})^2 + w_{0}^2 ) 
  \end{equation*}
   and
   \begin{equation*} 
   (u_{2}, v_{2}, w_{2}) =( |(2v_{0})^2 - w_{0}^2|,\quad 4 v_{0}w_{0}, \quad(2v_{0})^2 + w_{0}^2 )
   \end{equation*} 
   a three primitive Pythagorean triples satisfy the condition $ v_{0} v_{1} = u_{0}v_{2} $, then the parametrization of the edges
\begin{equation*}  
 a  =    u_{0}u_{2}= u_{0}(4v_{0}^2 -w_{0}^2) , \hspace{.5cm} b =  v_{0} u_{1}= v_{0}(4u_{0}^2-w_{0}^2), \hspace{.5cm} c = v_{0} v_{1} = 4 u_{0}v_{0}w_{0} 
\end{equation*}
give face diagonals
  \begin{equation*} 
{\displaystyle d= \sqrt{(u_{0}u_{2})^2+(v_{0}u_{1})^2 } = \sqrt{(u_{0}(4v_{0}^2 -w_{0}^2))^2+(v_{0}(4u_{0}^2-w_{0}^2))^2 } = w_{0}^3,   }
\end{equation*}
\begin{equation*} 
{\displaystyle   e = u_{0}w_{2} = u_{0}(4v_{0}^2+w_{0}^2) ,\hspace{2cm} f = v_{0}w_{1} = v_{0}(4u_{0}^2+w_{0}^2) }
\end{equation*}
leads to an Euler brick.
 \end{proof}

 \begin{footnotesize}
 \begin{table}[h!]
  \begin{tabular}{|c||c|c|c||c|c|c||c|c|c||c|c|c||c|c|c| } 
\hline
$ \# $ & $ u_{0} $ & $ v_{0} $ & $ w_{0} $   & $ u_{1} $ & $ v_{1} $ & $ w_{1} $  & $ u_{2} $ & $ v_{2} $ & $ w_{2} $ & $ a $ & $ b $ & $ c$  & $ d $ & $e $ & $ f$  \\ \hline
1  &  3  &  4 & 5 & 11 & 60 & 61 &  39 &  80 &  89 &  117 &  44 &  240 & 125 & 267 &  244     \\ \hline
2  & 5  & 12  & 13 & 69 & 260 & 269  & 407  & 624  & 745  & 2035  & 828  & 3120  & 2197  & 3725  &  3228   \\ \hline
3 & 15  & 8 & 17  & 611 & 1020  & 1189 & 33 & 544 & 545    & 495 & 4888  & 8160  & 4913  & 8175  & 9512 \\ \hline
4 & 7 & 24  & 25 &  429 &  700 & 821  & 1679 & 2400 & 2929   & 11753  & 10296   & 16800   & 15625 & 20503   & 19704 \\ \hline
5 & 21 & 20 & 29 &  923  & 2436  & 2605  & 759  & 2320  & 2441  & 15939  & 18460  &  48720  & 24389  &  51261 & 52100  \\ \hline
\end{tabular}
\bigskip
  \caption{The first five primitive Euler brick, given by the formula \eqref{eq13}, \eqref{eq14} and \eqref{eq17}, \eqref{eq18}. } 
 \end{table}
 \end{footnotesize}
 
 \begin{problem}
Let $ (u_{1}, v_{1}, w_{1}) $ and $ (u_{2}, v_{2}, w_{2}) $ two primitive Pythagorean triples and $ (u_{0}, v_{0}) $  two positive integers, related to Theorem \ref{th3}, If the condition $ (u_{0}u_{2})^2+(v_{0}u_{1})^2=\square $ is met, then $ v_{0}v_{1} = u_{0} v_{2} $ ?    
\end{problem}
 
 \section*{Conclusion}

 Let $ (u_{1},v_{1},w_{1}), (u_{2},v_{2},w_{2}), (u_{3},v_{3},w_{3}) $ be three primitive Pythagorean triples, and $ k_{1},k_{2},k_{3} $ are the $ \gcd $ of $ (a,b), (a,c), (b, c) $ respectively, in summary of the three theorems above, in the following table
  
 \begin{table}[h!]
  \begin{tabular}{|c| c|c|c||c |c|c|c||c|c|c|c|   } 
\hline
  & $ u_{1} $ & $ v_{1} $ & $ w_{1} $ &   & $ u_{2} $ & $ v_{2} $ & $ w_{2} $ &     & $ u_{3} $ & $ v_{3} $ & $ w_{3}$     \\ \hline
 $ \times k_{1} $  & $ k_{1} u_{1} $  &  $ k_{1} v_{1} $ & $ k_{1} w_{1} $   &   $ \times k_{2} $  & $ k_{2} u_{2} $ &   $ k_{2} v_{2} $ & $ k_{2} w_{2} $ & $ \times k_{3} $  & $ k_{3} u_{3} $   &  $ k_{3} v_{3} $ &    $ k_{3} w_{3} $        \\ \hline
  $ = $ &  $ a $ &  $ b $ & $ d $  &  =  &  $ a $ & $ c $  & $ e $  &  = &  $ b $ & $ c $ & $ f $     \\ \hline
    
\end{tabular}
\bigskip
  \caption{Summary of the three theorems } 
 \end{table}

\begin{problem}
Discover all the solutions for the Diophantine equation system 
 \begin{equation} \label{eq21} 
  (xu)^2+(yv)^2=\square, \hspace{3cm} (xv)^2+(yu)^2=\square,
 \end{equation}
   where $ (x,y,z) $ and $ (u,v,w) $ two primitive Pythagorean triples.
\end{problem}

\begin{conjecture}
we conjecture if the system \eqref{eq21} has a solution, then $ (x,y,z) $ or $ (u,v,w) $ equals $ (1,0,1) $.
\end{conjecture}
The guesswork was examined manually for triangles below 100.

\begin{example} $  $
 \begin{itemize}
\item If take  $ (x_{1},y_{1},z_{1})=(3,4,5) $ and $ (u_{1},v_{1},w_{1})=(5,12,13) $ , then 
\begin{equation*}   
  (x_{1}u_{1})^2+(y_{1}v_{1})^2= 2529 \neq \square, \hspace{1cm} (x_{1}v_{1})^2+(y_{1}u_{1})^2=1696 \neq \square.
 \end{equation*}
\item If take  $ (x_{2},y_{2},z_{2})=(20,21,29) $ and $ (u_{2},v_{2},w_{2})=(48,55,73) $ , then 
\begin{equation*}   
  (x_{2}u_{2})^2+(y_{2}v_{2})^2= 2255625  \neq \square, \hspace{1cm} (x_{2}v_{2})^2+(y_{2}u_{2})^2=1492^2 = \square.
 \end{equation*}
 
\item If take  $ (x_{3},y_{3},z_{3})=(11,60,61) $ and $ (u_{3},v_{3},w_{3})=(39,80,89) $ , then 
\begin{equation*}   
  (x_{3}u_{3})^2+(y_{3}v_{3})^2= 23224041 \neq \square, \hspace{1cm} (x_{3}v_{3})^2+(y_{3}u_{3})^2=2500^2 = \square.
 \end{equation*}

\item If take  $ (x_{4},y_{4},z_{4})=(7,24,25) $ and $ (u_{4},v_{4},w_{4})=(99,20,101) $ , then 
\begin{equation*}   
  (x_{4}u_{4})^2+(y_{4}v_{4})^2= 843^2 = \square, \hspace{1cm} (x_{4}v_{4})^2+(y_{4}u_{4})^2=5664976 \neq \square.
 \end{equation*}
 
\item  If take  {\small $ (x_{5},y_{5},z_{5})=(11,60,61)$} and {\small $(u_{5},v_{5},w_{5})=(17,144,145)$}, then 
\begin{equation*}   
  (x_{5}u_{5})^2+(y_{5}v_{5})^2= 74684569  \neq \square, \hspace{1cm} (x_{5}v_{5})^2+(y_{5}u_{5})^2=1884^2 = \square.
 \end{equation*}

\end{itemize}
\end{example}

Similarly, we make the following conjecture
\begin{conjecture} \label{cj2} 
   Let $ (x,y,z) $ and $ (u,v,w) $ two primitive Pythagorean triples. If the system of two Diophantine equations 
 \begin{equation} \label{eq22} 
  (zu)^2+(yv)^2=\square, \hspace{3cm} (zv)^2+(yu)^2=\square,
 \end{equation}
  has a solution, then $ (x,y,z) $ or $ (u,v,w) $ equals $(1,0,1)$.
\end{conjecture}

Any counterexample to conjecture \ref{cj2}, allows us to obtain the perfect Euler cuboid, edges

 \begin{equation}  
 a  =  xu , \hspace{2cm} b = xv, \hspace{2cm} c = yw \hspace{1cm}
\end{equation}
give face diagonals
\begin{equation}   
{\displaystyle d = xw,\hspace{1cm} e = \sqrt{(zu)^2+(yv)^2} ,\hspace{1cm} f = \sqrt{(zv )^2+(yu)^2 }  }
\end{equation}
and space diagonal  
\begin{equation}  
 g  = w z
\end{equation}


\bigskip


\begin{thebibliography}{1}  

\bibitem {1} Knill, Oliver, \href{https://people.math.harvard.edu/~knill/various/eulercuboid/lecture.pdf}{\it Treasure Hunting Perfect Euler bricks}, (February 24, 2009), Math table, Harvard University.
\bibitem {2} Knill, Oliver, \href{https://doi.org/10.48550/arXiv.2205.13285}{\it The Babylonian graph},   (May 26, 2022), arXiv:2205.13285v1 [math.CO].

\bibitem {3} Ian stewart, \href{ }{\it Visions of Infinity: The Great Mathematical Problems By Ian Stewart},  Chapter 17, Basic Books, New York, 2013.

\bibitem {4} Waclaw Sierpinski, \href{}{\it Pythagorean Triangles},  Yeshiva University, New York, 1962.

 \bibitem {5} Matson, Robert D, \href{http://unsolvedproblems.org/S58.pdf}{\it Results of a Computer Search for a Perfect Cuboid}.
\end{thebibliography}
\end{document}